\documentclass[11pt]{amsart}
\usepackage{amssymb,amsmath,amsfonts,amscd,euscript}

\newcommand{\nc}{\newcommand}

\nc{\g}{\mathfrak{g}}
\nc{\la}{\lambda}
\nc{\slth}{\widehat{\slt}}
\nc{\C}{\mathbb C }
\nc{\Pp}{\mathbb P }
\nc{\R}{\mathbb R }
\nc{\Z}{\mathbb Z }
\nc{\N}{\mathbb N }
\nc{\E}{\mathbb E }
\nc{\al}{\alpha}
\nc{\be}{\beta}
\nc{\ga}{\gamma}
\nc{\bp}{{\bf p}}
\nc{\bs}{{\bf s}}
\nc{\ve}{\varepsilon}
\nc{\ch}{{\mathop {\rm ch}}}
\nc{\Id}{{\mathop {\rm Id}}}
\nc{\ld}{\ldots}
\nc{\cd}{\cdots}
\nc{\hk}{\hookrightarrow}
\nc{\n}{\mathfrak{n}}
\nc{\T}{\otimes}

\nc{\pa}{\partial}
\nc{\bt}{{\bf t}}
\nc{\bq}{{\bf q}}
\nc{\Lc}{\EuScript L}
\nc{\Hz}{\EuScript H}
\nc{\Fh}{\EuScript F}
\nc{\Gh}{\EuScript G}
\nc{\Mg}{\EuScript M}
\nc{\Ll}{\EuScript L}

\nc{\at}[2]{\genfrac{}{}{0pt}{}{#1}{#2}}

\nc{\beq}{\begin{equation}}
\nc{\eeq}{\end{equation}}

\newtheorem{theo}{Theorem}[section]
\newtheorem*{theo*}{Theorem}
\newtheorem{lem}{Lemma}[section]

\newtheorem{cor}{Corollary}[section]
\newtheorem{rem}{Remark}[section]

\begin{document}
\author{E.Feigin}
\title
[Formal genus $0$ theories]
{$N=1$ formal genus $0$ Gromov-Witten theories and Givental's formalism}

\address{Evgeny Feigin:\newline
{\it Mathematical Institute, Cologne University,
Weyertal 86-90, 50931, Cologne, Germany}
\newline and \newline
{\it Tamm Theory Division, Lebedev Physics Institute,
Russian Academy of Sciences,\newline
Russia, 119991, Moscow, Leninski prospect, 53}}
\email{evgfeig@gmail.com}

\begin{abstract}
In \cite{Gi3} Givental introduced and studied a space of formal genus zero Gromov-Witten theories $GW_0$,
i.e. functions satisfying string and dilaton equations and topological
recursion relations. A central role in the theory  plays the geometry of certain Lagrangian cones
and a twisted symplectic group of hidden symmetries.
In this note we show that the Lagrangian cones description of the action of this group coincides with
the genus zero part of Givental's quantum Hamiltonian formalism. As an application we identify
explicitly the space of $N=1$ formal genus zero GW theories with lower-triangular
twisted symplectic group modulo the string flow.
\end{abstract}

\maketitle

\section*{Introduction}
Let $X$ be a compact Kahler manifold and $H=H^*(X)$ be the space of cohomology
of $X$ equipped with a non-degenerate Poincare pairing $(\cdot,\cdot)$.
Let $\Hz_+=H\T\C[z]$ be the big phase space. The genus $g$ Gromow-Witten
descendant  potential of $X$ is a
generating function $F_g^X$ for Gromov-Witten invariants (\cite{LT}). The series $F_g^X$
depends on variables $\{t_i\}_{i=0}^\infty$, $t_i\in H$ and thus can be considered
as a function on the big phase space.

There exists a large family of relations which are satisfied by all descendant
potentials (\cite{BP}, \cite{EX}, \cite{FSZ}, \cite{Ge1}, \cite{Ge2}, \cite{L1},
\cite{LSS}, \cite{Ma}). These relations come from the geometry of the moduli spaces
$\overline{\Mg}_{g,n}$.
They include the string and the dilaton equations (see \cite{W1}), and the set of topological recursion
relations (see \cite{Ge1}).
In \cite{Gi2}, \cite{Gi3} axiomatic approach to the study of
the so called formal Gromow-Witten theories is developed.
Namely the idea is to capture the main properties of geometric Gromow-Witten
potentials and to consider the space of functions satisfying those properties.
Recent applications of Givental's formalism include Frobenius structures
(see \cite{L1}, \cite{L2}, \cite{L3}, \cite{LP}),  $r$-spin
Gromow-Witten theory (see \cite{FSZ},\cite{CZ}), integrable systems (\cite{Gi4},
\cite{FSZ}, \cite{Mi}), topological and cohomological field theories (\cite{T}, \cite{K}, \cite{S}).

The simplest and the most concrete is the genus zero case.
The space of formal genus zero Gromow-Witten theories
$GW_0$ is defined as follows.
Let $H$ be an $N$-dimensional vector space equipped with a non-degenerate
bilinear form $(\cdot,\cdot)$ and a distinguished vector ${\bf 1}\in H$.
Consider a space of functions
$F(t_0,t_1,\dots)$ satisfying the string equation, the dilaton equation
and the set of genus zero topological recursion relations.
In particular, genus zero parts of geometric Gromov-Witten
potentials and genus zero parts of Witten's $r$-spin potentials
(see \cite{W1}, \cite{JKV}) are elements of $GW_0$.
It is proved in \cite{CG}, \cite{Gi1} that the space of formal genus zero theories
can be equipped with an action of infinite-dimensional symplectic Lie
groups $G_+$ and $G_-$ (the so-called upper- and lower-triangular twisted loop groups).
The existence of these
groups of hidden symmetries  is a very powerful tool for the study
of formal theories (see \cite{LP}, \cite{L3} and references within).

There exist two definitions of the action of twisted loop groups.
The first one is based
on the quantization of certain Hamiltonians (see \cite{Gi2})
and goes through the infinitesimal action of the corresponding Lie algebras
on the space of total descendant potentials. We recall
some details in Section~$1$.
The second definition is based on
the Lagrangian cones formalism (see \cite{CG}, \cite{Gi3}). Namely let
$\Hz=H\T\C[z,z^{-1}]]$ be the space of Laurent series with values in $H$
 equipped with a symplectic form
$$\Omega(f(z),g(z))=\frac{1}{2\pi i} \oint (f(-z),g(z))dz.$$
The Lagrangian decomposition $\Hz=\Hz_+\oplus H\T z^{-1}\C[[z^{-1}]]$ ($\Hz_+$ is a big phase space)
leads to the identification of $\Hz$ with the cotangent bundle $T^*\Hz_+$.
For $F\in GW_0$ let $\Lc_F\hk \Hz$ be a graph of $dF$ (recall that $F$ can be viewed as a function
on $\Hz_+$). In \cite{Gi3}
the system $DE+SE+TRR$ is rewritten in terms of geometric properties of $\Lc_F$. This
approach allows to define the action of $G_\pm$ on $GW_0$ as well as to study
the properties of this action on the subspace $GW_0^{ss}$ of semi-simple theories.

In this paper we suggest another -- more algebraic -- approach to the study
of the space of formal genus $0$ Gromov-Witten theories.
The main points are the exact formula for the action of the lower-triangular twisted loop
group and the form of the expansion of elements of $GW_0$ as Laurent series
in a distinguished variable. We will deal only with the action of the lower-triangular
group $G_-$ since the "opposite" group $G_+$ doesn't preserve the space $\Hz$ (see \cite{Gi3}
for the discussion of possible completitions). We briefly describe our approach below.

First, starting
from the Lagrangian cones formalism we derive exact
formula for the action of $G_-$.
The resulting  formula coincide with a genus zero restriction of the Givental's
quantum action of the twisted loop group on the space of total descendent
potentials (see \cite{Gi2}). Thus we prove that the Lagrangian cones approach
and quantization procedure produce  the same actions of $G_-$
on genus $0$ formal theories (see also \cite{Gi3}, \cite{L1}, \cite{L3}).

Next, using the dilaton equation we write elements of $GW_0$ as Laurent series in
a distinguished variable.
Namely we fix a basis $\{\phi_\al\}_{\al=1}^N$ of $H$ such that
$\phi_1={\bf 1}$. This gives a basis $\phi_\al\T z^n$ of $\Hz_+$ and the corresponding
coordinates $t^\al_n$. We show that
any element $F\in GW_0$ can be written in a form
\beq
F=\sum_{n\ge 0} (t_1^1-1)^{2-n} c_n,
\eeq
where $c_n$ are homogeneous degree $n$ polynomials in variables $t^\al_n$,
$(\al,n)\ne (1,1)$.
As an application we study $N=1$ case. We show explicitly that the group
$G_-$ acts transitively and, modulo the string flow, freely on the space $GW_0$.
We note that for $N=1$ the group $G_+$ acts trivially on $GW_0$ (see
\cite{Gi2}, \cite{FP}). Thus our results agree with Givental's theorem
which states that for general $N$ one needs the action of both groups $G_\pm$
in order to generate the whole space $GW_0$ starting from the Gromov-Witten
potential of $N$ points.

We finish the introduction part with a following  remark. The results
of our paper seems to  be known to experts. The novelty is an algebraic
approach to the Givental's theory. The advantage of this approach is two-fold.
Firstly, it makes dome constructions  more clear and allows to prove
certain statements or simplify  the existed proofs. Secondly, we hope that
our approach can be applied to the study of non-semisimple GW theories, which
are not covered by the geometric theory (see a discussion in \cite{L3}).

Our paper is organized as follows.
In Section $1$ we recall main points of the Givental's formalism.
In Section $2$ we derive exact formulas for the action of the lower-triangular
subgroup on $GW_0$.
In Section $3$ we recall the connection with Frobenius structures.
In Section $4$ we describe explicitly the space $GW_0$ for $N=1$.

\noindent{{\bf Acknowledgements}.}
EF thanks M.Kazarian, S.Shadrin and D.Zvonkine for useful discussions on Givental's
theory.
This work was partially supported by the RFBR  Grants
06-01-00037, 07-02-00799, by Pierre Deligne fund based on his 2004
Balzan prize  in mathematics, by Euler foundation and by Alexander von Humboldt
Fellowship.

\section{Givental's formalism}
In this section we recall main points of the formalism developed in
\cite{CG}, \cite{Gi1},\cite{Gi2}.

Let $H$ be an $N$-dimensional vector space equipped with a
symmetric
non-degenerate bilinear form $(\cdot,\cdot)$ and with a distinguished
nonzero element ${\bf 1}$.
Let $\Hz$ be the space
of Laurent series $H((z^{-1}))$ in $z^{-1}$. This space carries a
symplectic form $\Omega$:
$$\Omega(f(z),g(z))=\frac{1}{2\pi i} \oint (f(-z),g(z))dz.$$
The space $\Hz$ is naturally a $\C[z,z^{-1}]$-module. Any operator
$M(z)\in\mathrm{End}(\Hz)$, which commutes with the action of $z$
can be written in a form $M(z)=\sum_{i\in\Z} M_iz^i$. If $M(z)$
preserves $\Omega$ then one has
$$M^*(-z)M(z)=\Id,$$ where $\phantom{}^*$ is  adjoint
with respect to $(\cdot,\cdot)$.
An algebra of $\mathrm{End}(H)$-valued Laurent series $\mathrm{End}(H)((z^{-1})$
acts on the space $\Hz$.

\begin{lem}\label{L}
Let $M(z)\in \mathrm{End}H((z^{-1})$ be a series such that
$M^*(-z)M(z)=\Id$. Then $M_i=0$ for $i>0$.
\end{lem}
\begin{proof}
Suppose $M_i=0$ for $i>N>1$. Then because of $M^*(-z)M(z)=\Id$  one gets
$M_N^* M_N=0$. Therefore $M_N=0$ since $(\cdot,\cdot)$ is non-degenerate.
\end{proof}

The group of operators $S(z)=\Id+\sum_{i> 0} S_iz^i$ satisfying
$S^*(-z)S(z)=\Id$ is called {\it the lower-triangular group}. We denote it
by $G_-$. The elements of this group are sometimes called the calibrations of
Frobenius manifolds (see Section $3$). They are also involved
in the ancestor-descendent potentials correspondence (see \cite{Gi2},\cite{KM}).

\begin{rem}
The upper-triangular group $G_+$ consists of the operators
$$R(z)=\Id+\sum_{i>0} R_iz^i.$$
These operators play the central role in the Givental's
theory of Frobenius manifolds and formal
Gromov-Witten potentials (see \cite{Gi1}).
We are not considering the group $G_+$ in this paper because of
two reasons. First because of Lemma $\ref{L}$ the upper-triangular group  does not
exactly fit the Lagrangian cone formalism (see \cite{CG}, \cite{Gi3} for
discussion of possible modifications). Another reason is that in our main application ($N=1$)
we only need the action of the group $G_-$ (because of  Faber-Pandharipande relations,
see \cite{Gi2},\cite{FP}).
\end{rem}

We now describe the action of the Giventals group on the axiomatic genus $0$
theories and the quantization formalism.
Let $\{\phi_\al\}$ be a basis of $H$, such that
$\phi_1={\bf 1}$, $(\phi_\al,\phi_\beta)=g_{\al,\beta}$.
Let $t_n$ be an element of $H\T z^n$. Then we set
$t_n=\sum_\al t_n^\al \phi_\al$. In what follows we deal with functions
$F(t_0,t_1,\dots)$, $t_i\in H$. We write $\pa_{\al,n}$ for the partial
derivative  with respect to $t_n^\al$.
As usual, $(g^{\al,\be})=(g_{\al,\be})^{-1}$.

By definition, a function
$F(t_0,t_1,\dots)$, $t_i\in H$ is called axiomatic genus zero Gromov-Witten
theory if it satisfies the string equation (SE)
$$
\sum_{i\ge 0} t^\al_i\pa_{\al,i} F(\bt)  - \pa_{1,1} F(\bt)= 2F(\bt),
$$
the dilaton equation (DE)
$$
\sum_{i\ge 0} t^\al_i\pa_{\al,i} F(\bt)  - \pa_{1,1} F(\bt)= 2F(\bt),
$$
and the system
of topological recursion relations (TRR), labeled by triples of
integers $k,l,m\ge 0$ and $\al,\be,\gamma=1,\dots,\dim H$:
$$
\pa_{\al,k+1}\pa_{\be,l}\pa_{\ga,m} F(\bt)=
\pa_{\al,k}\pa_{\mu,0} F(\bt) g^{\mu,\nu}
\pa_{\nu,0}\pa_{\be,l}\pa_{\ga,m} F(\bt).
$$
We denote the set of functions satisfying string, dilaton and
topological recursion relations by $GW_0$.
It is convenient to introduce the dilaton shifted variables
$$q^\al_n=t^\al_n-\delta_{\al,1}\delta_{n,1}.$$
In new variables the system $(DE)$, $(SE)$ and $(TRR)$ for the function
$F(\bq)=F(q_0,q_1,\dots)$ reads as
$$
\sum_{i\ge 0} q^\al_i\pa_{\al,i} F(\bq) = 2F(\bq),\eqno (DE)
$$
$$
\sum_{i\ge 0} q^\al_{i+1}\pa_{\al,i} F(\bq) =-\frac{1}{2}(q_0,q_0), \eqno (SE),
$$
$$
\pa_{\al,k+1}\pa_{\be,l}\pa_{\ga,m} F(\bq)=
\pa_{\al,k}\pa_{\mu,0} F(\bq) g^{\mu,\nu}
\pa_{\nu,0}\pa_{\be,l}\pa_{\ga,m} F(\bq).
\eqno (TRR)
$$

\begin{lem}   \label{q_1}
Any axiomatic genus zero theory $F$ can be written in a form
\begin{equation}
\label{2-n}
F(\bq)=\sum_{n\ge 0} (q^1_1)^{2-n} c_n,
\end{equation}
where $c_n$ are degree $n$ series in variables $q^\al_n$,
$(\al,n)\ne (1,1)$.
\end{lem}
\begin{proof}
Let $$t^{\bf \al}_{\bf k}=t^{\al_1}_{k_1}\dots t^{\al_n}_{k_n}$$
be a monomial satisfying
$(\al_i,k_i)\ne (1,1)$. Let $a^{\bf \al}_{\bf k}(l)$ be a coefficient
of  $t^{\bf \al}_{\bf k} (t^1_1)^l$ in $F(t_0,t_1,\dots)$.
Comparing the coefficients of $t^{\bf \al}_{\bf k} (t^1_1)^l$ in the
right and left hand sides of the dilaton equation one gets
$$a^{\bf \al}_{\bf k}(l+1)=a^{\bf \al}_{\bf k}(l)\frac{l+n-2}{l+1}.$$
But the same relation holds for the coefficients $b(l)$ of $t^l$ in the
expansion of $(1-t)^{2-k}$. Now Lemma follows from the relation
$q^1_1=t^1_1-1$.
\end{proof}

\begin{rem}
For the genus $0$ Gromov-Witten potential of the point $F^{pt}_0$ we have
$$c_3=-\frac{1}{6} (q^1_0)^3,\ c_0=c_1=c_2=0,$$
because the Deligne-Mumford spaces
$\overline{\Mg}_{0,0}$, $\overline{\Mg}_{0,1}$ and  $\overline{\Mg}_{0,2}$
are missing.
In general, the coefficients $c_0$, $c_1$ and $c_2$ do not vanish.
\end{rem}

In \cite{Gi3} Givental constructed an action of the group of $\Omega$-preserving
operators on the space $GW_0$. We briefly outline his construction.
Let $n$ be some nonnegative integer. Recall
the coordinates $q_n^\al$ in the space $H\T z^n$, $n\ge 0$. We denote by $p_n^\al$ the Darboux
coordinates with respect to $\Omega$
in the space $H\T z^{-n-1}$ and set
$$
q_n=\sum_{\al=1}^N q_n^\al \phi_\al,\
p_n=\sum_{\al=1}^N p_n^\al \phi^\al,
$$
where $\phi^\al=\sum_{\be=1}^N g^{\al,\be} \phi_\be$.
Thus any element of $\Hz$ can be written in a form
$$\sum_{n\ge 0} q_n z^n  + \sum_{n\ge 0} p_n (-z)^{-n-1}.$$
For any $F\in GW_0$ let $\Lc_F\subset \Hz$ be a graph of $dF$. Namely,
\beq
\label{Lc}
\Lc_F=\left\{\sum_{i\ge 0} p_i(-z)^{-i-1}+  \sum_{i\ge 0} q_i z^i:\
p_i^\al=\frac{\pa F}{\pa q_i^\al}\right\}.
\eeq
Because of the dilaton equation $(DE)$ the subset $\Lc_F$ is a cone.

\begin{theo}\cite{Gi3}
\label{G}
A cone $\Lc\subset \Hz$ is equal to some $\Lc_F$ for $F\in GW_0$
if and only if $\Lc$ is Lagrangian and for any $f\in\Lc$ the
tangent space $L_f=T_f\Lc$ is tangent to $\Lc$ exactly along $zL$.
\end{theo}

\begin{cor}
The group $G_-$ acts on $GW_0$.
\end{cor}

The following simple Lemma will be important for our computations:
\begin{lem}  \label{rec}
Let $F$ be an element of $GW_0$ and let $\Lc_F\hk \Hz$ be the corresponding
cone with the projection $\pi:\Lc_F\to\Hz_+$. Then
$$2F(\bq)=\sum_{i\ge 0} (p_i,q_i),\ \bp=\pi^{-1}\bq.$$
\end{lem}
\begin{proof}
Follows from the definition of $\Lc_F$ and the dilaton equation.
\end{proof}

Another description of the action of $G_-$ can be given in terms of the
genus zero restriction of the Givental's quantization formalism.
This theory defines an action of the groups $G_\pm$ on the space of
total descendent potentials, involving all genera.
Namely, let $\g_+$, $\g_-$  be Lie
algebras of the groups $G_+$, $G_-$.
Explicitly,
$$\g_\pm=\{\sum_{i>0} a_iz^{\pm i}: a_i^*=(-1)^{i+1} a_i\}.$$
For any $a\in\g_\pm$ consider a Hamiltonian
$$H_a(f)=\frac{1}{2}\Omega(af,f).$$
The element $f\in\Hz$ can be written in a form
$f=\sum_{i\ge 0} p_i (-z)^{-i-1}+\sum_{i\ge 0} q_i z^i$, $p_i,q_i\in H$.
Thus $H_a$ is a quadratic function in variables
$p_i^\al$, $q_j^\be$.
The quantization rule
$$p_i^\al p_j^\be\mapsto \hbar \partial_{i,\al}\partial_{j,\be},\
p_i^\al q_j^\be\mapsto  q_{j,\be}\partial_{i,\al},\
q_i^\al q_j^\be\mapsto \hbar^{-1} q_{i,\al}q_{j,\be}$$
defines operators $\widehat H_a$ on $\hbar$-series
with values in the space of functions on $\Hz_+=\Hz\T\C[[z]]$.
The exact formulas are given in \cite{L1} (see also \cite{FSZ}).
Namely
\begin{gather}
\label{Sa}
\widehat{H_{s_lz^{-l}}}=
\sum_{n\ge 0}\sum_{\al,\be}
(s_l)_{\al,\be} q_{l+n}^\al \pa_{n,\be} +
\frac{1}{\hbar} \sum_{n=0}^{l-1} \sum_{\al,\be} (s_l)_{\al,\be}
(-1)^n q_n^{\al}q_{l-1-n}^\be,\\
\label{Ra}
\widehat{H_{r_lz^l}}=
\sum_{n\ge 0}\sum_{\al,\be}
(r_l)_{\al,\be} q_n^\al \pa_{l+n,\be} +
\frac{\hbar}{2} \sum_{n=0}^{l-1} \sum_{\al,\be} (r_l)_{\al,\be}
(-1)^{n+1}\pa_{n,\al}\pa_{l-1-n,\be},
\end{gather}
where $l>0$ and $(a)_{\al,\be}$ are entries of a matrix $a$ in the basis
$\phi_\al$ of $H$.
These formulas give rise to the action of the groups $G_+$ and $G_-$ on the
space of formal total descendant potentials of the form
$$\Fh=\exp(h^{-1}F_0+F_1+hF_2+\dots).$$.

In the following Corollary we extract the genus zero part of the action of
$G_-$. We use the notation $\bq(z)=\sum_{i\ge 0} q_iz^i$ and write $F(\bq(z))$
instead of $F(\bq)$.

\begin{cor} \label{Sq}
Let $F$ be a genus $0$ part of a total descendent potential $Z$. Then the genus zero
part of the image $\widehat S(z) \Fh$ is given by the formula
\beq
(\widehat S(z)F)(\bq(z))=\frac{1}{2}\sum_{k,l\ge 0} (W_{k,l}q_k, q_l)+
F([S^{-1}(z)\bq(z)]_+),\\
\eeq
where the operators $W_{k,l}$ are defined by
\beq\label{W1}
\sum_{k,l\ge 0} \frac{W_{k,l}}{z^kw^l}=\frac{(S^*)^{-1}(w)S^{-1}(z)-\Id}{z+w}
\eeq
and $[S^{-1}(z)\bq(z)]_+$ is a power series truncation of the  series
$S^{-1}(z)\bq(z)$.
\end{cor}
\begin{proof}
Follows either from the formula $(\ref{Sa})$ or from the genus zero restriction
of the formula from \cite{Gi2}, Proposition $5.3$.
\end{proof}

\begin{rem}
The change of variables $\bq(z)\to [S^{-1}(z)\bq(z)]_+$ is not always well-defined
operation. One should be careful to avoid infinite coefficients in the decomposition
$(\ref{2-n})$ after the substitution.
In the following Lemma we formulate a sufficient condition.
\end{rem}

\begin{lem}\label{S1}
Let $S_1=0$. Then the expression $F_0([S^{-1}(z)\bq(z)]_+)$ is well-defined.
\end{lem}
\begin{proof}
We consider a summand $(q_1^1)^{2-n}c_n$ from the decomposition $(\ref{2-n})$, where
$c_n$ is a degree $n$ polynomial in variables $q_n^\al$, $(\al,n)\ne (1,1)$.
After the change of variables ${\bf q}(z)\to [S^{-1}(z)\bq(z)]_+$ the series   $c_n$ remains
independent of $q_1^1$ because of the condition $S_1=0$. Therefore the sum
$$(\sum_{n\ge 0} (q_1^1)^{2-n} c_n)([S^{-1}(z)\bq(z)]_+)$$
is well-defined.
\end{proof}

\section{"Negative" (lower triangular) subgroup}
\label{S(z)}
Let $S(z)=1+S_1z^{-1}+S_2z^{-2}+\dots$ be an element of the
group $G_-$,
\beq
\label{sympl}
S^*(-z)S(z)=1.
\eeq
For instance
$$S_1^*+S_1=0,\ S_2^* + S_1^* S_1 +  S_2=0.$$
Let $F(\bq)$ be a formal theory, $\Lc_F$ be the corresponding  Lagrangian
cone.
We use Lemma  $\ref{rec}$ to compute the action of $G_-$ on $GW_0$.
By definition
$S(z)(\bp,\bq)=(\bp',\bq')$, where
\begin{gather}
q_i'= \sum_{n\ge 0} S_n q_{i+n}, \label{q'}\\
p_i'=\sum_{n=0}^i (-1)^n S_n p_{i-n} - (-1)^i \sum_{n\ge 0} S_{i+n+1} q_n
\label{p'}.
\end{gather}

\begin{lem}
\beq
\label{p'q'}
\sum_{i\ge 0} (p_i', q_i')= \sum_{i\ge 0} (p_i, q_i) -
\sum_{i\ge 0} (-1)^i (S_{i+1}q_0+S_{i+2}q_1+\dots, q_i').
\eeq
\end{lem}
\begin{proof}
Using the formulas $(\ref{q'})$ and $(\ref{p'})$ we obtain
$$\sum_{i\ge 0} (p_i', q_i')=
\sum_{i\ge 0}
\left(\sum_{n=0}^i (-1)^n S_n p_{i-n} +(-1)^{n+1}\sum_{n\ge 0} S_{i+n+1} q_n,
q_i'\right).$$
The direct computation shows that because of the relation
$S^*(-z)S(z)=1$ all the terms  of the form $(p'_i,q'_j)$ in the right hand
side vanish.  The rest terms are exactly the right hand side of
$(\ref{p'q'})$.
\end{proof}

\begin{lem}
\label{W}
Let $W_{k,l}$, $k,l\ge 0$ be operators defined by
\beq
\label{decomp}
\frac{S^{-1*}(w) S^{-1}(z)-1}{z+w}=
\sum_{k,l\ge 0} \frac{(W_{k,l} q_k,q_l)}{z^kw^l}.
\eeq
Then for any $i\ge 0$
$$
(-1)^{i+1}\sum_{n\ge 0} S_{i+n+1} q_n=
\sum_{k\ge 0} W_{k,i} q_k'.$$
\end{lem}
\begin{proof}
Using the definition of $q_k'$ we obtain
\beq
\sum_{k\ge 0} W_{k,i} q_k'=\sum_{k\ge 0} (\sum_{m=0}^k W_{m,i} S_{k-m}) q_k.
\eeq
To compute  $\sum_{m=0}^k W_{m,i} S_{k-m}$  we multiply
$(\ref{decomp})$ by $S(z)$ from the right:
\beq
\label{Wkl}
\sum_{k,l\ge 0} \frac{W_{k,l}}{z^kw^l} S(z)=
\frac{S^{-1*}(w)S^{-1}(z)-1}{w+z}S(z)=\frac{S(-w)-S(z)}{w+z}.
\eeq
We now decompose the rightmost and leftmost expressions  in $(\ref{Wkl})$.
\begin{gather*}
\sum_{k,l\ge 0} \frac{W_{k,l}}{z^kw^l} S(z)=
\sum_{k,l\ge 0} \frac{1}{z^k w^l} \sum_{m=0}^k W_{m,l} S_{k-m}, \\
\frac{S(-w)-S(z)}{w+z}=
\sum_{k,l\ge 0} \frac{1}{z^k w^l} (-1)^{l+1} S_{k+l+1}.
\end{gather*}
We conclude that  $\sum_{m=0}^k W_{m,l} S_{k-m}=(-1)^{l+1} S_{k+l+1}.$
Lemma is proved.
\end{proof}

\begin{theo}
\label{SLagr}
$a)$. \ The action of the group $G_-$ is given by
$$(\widehat S(z)F)(\bq)=F([S^{-1}(z)\bq(z)]_+)+
\frac{1}{2}\sum_{k,l\ge 0} (W_{k,l}q_k,q_l),$$
where  $[S^{-1}(z)\bq(z)]_+$ is a power series truncation.\\
$b).$\ The action of the Lie algebra $\g_-$  is given by the genus zero restriction of
the formula $(\ref{Sa})$.
\end{theo}
\begin{proof}
Formula $(\ref{p'q'})$ and Lemma $\ref{W}$ show that
$$(S(z) F)(\bq')=F(\bq(z))+ \frac{1}{2} W(\bq',\bq').$$
Now part $a)$ of our theorem follows from the equation
$$\sum_{i\ge 0} q_i'z^i=[S(z)\bq(z)]_+.$$
The statement $b)$ of the theorem follows from the part $a)$.
\end{proof}

\begin{cor}
The definition of the action of $G_-$ via the Lagrangian cones formalism
is equivalent to those given via the genus zero restriction  of the
Givental's quantization procedure.
\end{cor}
\begin{proof}
Follows from  Corollary  $\ref{Sq}$ and
Theorem $\ref{SLagr}$.
\end{proof}

\section{Frobenius structures.}
In \cite{Gi3} Givental used the Lagrangian cones technique to
study Frobenius structures on $H$ (\cite{D}). In particular he
proved that each element of $GW_0$ determines a Frobenius manifold
structure on $H$. In this section we reprove this statement
algebraically.
\begin{lem}
Let $F\in GW_0$. Then the
function $\Phi(q)$, $q\in H$ defined by
$$\Phi(q)=F(q,0,0,\dots)$$
satisfies the WDVV equation.
\end{lem}
\begin{proof}
Let $1\le a,b,c,d\le N=\dim H$. We introduce the notation
\beq
D=\pa_{a,0}\pa_{b,1}\pa_{c,0}\pa_{d,0} F.
\eeq
We denote by $TRR_{\al,k;\be,l;\gamma,m}$ the topological recursion
relation which corresponds to the triples $(k,l,m)$, $(\al,\be,\gamma)$.
Differentiating both sides of the relation $TRR_{b,1;c,0;d,0}$ with respect to $q_0^a$
we obtain
\beq
\label{bad}
D=\pa_{a,0}\pa_{b,0}\pa_{e,0} F g^{ef} \pa_{f,0}\pa_{c,0}\pa_{d,0} F+
\pa_{b,0}\pa_{e,0} F g^{ef} \pa_{f,0}\pa_{a,0}\pa_{c,0}\pa_{d,0} F.
\eeq
Similarly, differentiating both sides of
$TRR_{b,1;a,0;d,0}$ with respect to $q_0^c$ we obtain
\beq
\label{bac}
D=\pa_{c,0}\pa_{b,0}\pa_{e,0} F g^{ef} \pa_{f,0}\pa_{a,0}\pa_{d,0} F+
\pa_{b,0}\pa_{e,0} F g^{ef} \pa_{f,0}\pa_{a,0}\pa_{c,0}\pa_{d,0} F.
\eeq
The equality of the right hand sides of $(\ref{bad})$ and $(\ref{bac})$ imply
the WDVV equation for $\Phi(q)$:
$$\pa_a\pa_b\pa_e \Phi g^{ef} \pa_f\pa_c\pa_d \Phi=
  \pa_c\pa_b\pa_e \Phi g^{ef} \pa_f\pa_a\pa_d \Phi.$$
\end{proof}

Lemma above is equivalent to the following Corollary
(see \cite{Gi3}).
\begin{cor}
Each $F\in GW_0$ determines a Frobenius manifold structure on $H$ via
the multiplication
$$\phi_{\al}\bullet_{t_0}\phi_\be= \sum A_{\al\be}^\gamma(F) \phi_{\gamma},\
A_{\al,\be}^\gamma(F)=\sum_{\la=1}^N
g^{\gamma\la}(\pa_{\la,0}\pa_{\al,0}\pa_{\be,0} F)_{t_1=t_2=\dots=0}.$$
\end{cor}

The exact formulas for the action of $G_-$ provide an immediate corollary
(see \cite{Gi3}).

\begin{cor}
Let $S(z)$ be an element of $G_-$ such that $S_1=0$ (see Lemma $\ref{S1}$).
Then the action of $S(z)$ does not change the Frobenius structure.
\end{cor}
\begin{proof}
By definition
$$A_{\al\be}^\gamma(\widehat S(z)F)=
\sum_{\la=1}^N g^{\gamma\la} \pa_{0,\la}\pa_{0,\al}\pa_{0,\be}
(F([S^{-1}(z)\bq(z)]_+)),$$
where derivatives are taken after the substitution $t_1=t_2=\dots=0$.
But
$$\left(F([S^{-1}(z)\bq(z)]_+)-F(\bq(z))\right)_{t_1=t_2=\dots=0}=0.$$
Corollary is proved.
\end{proof}

\section{The rank $1$ case}
In this section we consider the case $N=\dim H=1$.
We show that the space $GW_0$ form a single orbit of the lower-triangular
group.

We first recall that the genus $0$  Gromov-Witten potential of a point is defined by
$$
F^{pt}_g(t_0,t_1,\dots)=\sum_{n\ge 3} \frac{1}{n!} \sum_{i_1,i_2,\dots,i_n\ge 0} t_{i_1}\dots t_{i_n}
\int_{\overline{\Mg}_{0,n}} \psi_1^{i_1}\dots \psi_n^{i_n},
$$
where $\psi_i$ denotes the first  Chern class of the line bundle on
$\overline{\Mg}_{0,n}$, whose fiber at the point $(C, x_1,\dots, x_n)$
is equal to $T^*_{x_i} C$.
The Laurent series expansion of $F^{pt}_0(q_0,q_1,\dots)$ in variable
$q_1$ starts with
$$F_0^{pt}=-\frac{q_0^3}{6q_1}-\frac{1}{24}\frac{q_0^4q_2}{q_1^3}+\dots.$$

It is easy to see that if $N=1$ then
$$\g_-=\{\sum_{i\ge 1} a_iz^{-2i+1},\ a_i\in\C\},\
G_-=\{\exp(\sum_{i\ge 1} a_iz^{-2i+1}),\ a_i\in\C\}$$
are abelian Lie algebra and Lie group.
The topological recursion
relations are labeled by triples of nonnegative integers. We set $g_{1,1}=1$
and denote the relation
$$
\pa_{k+1}\pa_l\pa_m F(\bq)=
\pa_k\pa_0 F(\bq) \pa_0\pa_l\pa_m F(\bq)
$$
by $TRR_{k,l,m}$. We also use a following notation:
let $(A)$ be an equality of two series in the variable $q_1$:
$$\sum_{n\in\Z} G_nq_1^{n}= \sum_{n\in\Z} H_nq_1^{n}. \eqno (A)$$
We denote an equality $G_n=H_n$ by $(A(n))$.

Recall (see Lemma $\ref{q_1}$) that any element $F\in GW_0$
can be written in a form
\begin{equation}
F(\bq)=\sum_{n\ge 0} q_1^{2-n} c_n(q_0,q_2,q_3,\dots),
\end{equation}
where $c_n$ are degree $n$ series in variables $q_0$, $q_2$, $\dots$.
Consider the string equation $(SE)$.
We note that equations $(SE(n))$ are trivial for
$n>2$  and the first $3$ nontrivial equations  are of the
form:
$$\pa_0 c_1=0, \eqno (SE(2)),$$
$$\pa_0 c_2+\sum_{i\ge 2} q_{i+1}\pa_i c_1=0, \eqno (SE(1)),$$
$$\pa_0 c_3+q_2 c_1+\sum_{i\ge 2} q_{i+1}\pa_i c_2=-\frac{1}{2}q_0^2,
\eqno (SE(0)).$$
In general  the following Lemma holds.
\begin{lem}
An explicit form of $(SE(n))$, $n<0$ is given by
$$\pa_0 c_{-n+3}+(-n+1)q_2 c_{-n+1}+\sum_{i\ge 2} q_{i+1}\pa_i c_{-n+2}=0.$$
\end{lem}

It turns out that among the TRR equations we only need the $TRR_{0,1,1}$.
In the following lemma we write this equation explicitly.
\begin{lem}
Equations $TRR_{0,1,1}(n)$ are trivial for $n>-4$.
Equations $TRR_{0,1,1}(n)$, $n\le -4$ are given by
$$(n+3)(n+2)(n+1)c_{-n-1}=\sum_{n_1+n_2=2-n} (2-n_2)(1-n_2)
\pa_0^2 c_{n_1} \pa_0 c_{n_2}.
$$
\end{lem}

\begin{lem}
For any $F\in GW_0$ the coefficient of $q_1^{-1}$ is a cube of some linear form.
That is
\beq
c_3=-\frac{1}{6}(q_0+\al_2q_2+\al_3q_3+\dots)^3
\eeq
for some constants $\al_i$, $i\ge 2$.
\end{lem}
\begin{proof}
Using relations $(SE(0))$ and $(SE(1))$ we obtain $\pa_0^2 c_2=0$.
Therefore the relation $(TRR_{0,1,1}(-4))$ reads as
$$
-6c_3=2\pa_0^2 c_3\pa_0 c_3.
$$
Using the relation $\pa_0^3c_3=1$ (which comes from $(SE(0))$)
we derive that
$$6c_3=(\pa_0^2 c_3)^3.$$
Lemma follows.
\end{proof}

\begin{cor}
The group $G_-'=G_-/\{\exp(az^{-1}), a\in\C\}$ acts freely on the set $GW_0$.
\end{cor}
\begin{proof}
We first recall that because of the string equation the action of  $z^{-1}\in \g_-$
on $GW_0$ is trivial. Next, Theorem $\ref{SLagr}$ shows that the coefficient of
$q_1^{-1}$ in $\exp(-az^{2k+1})F$, $k\ge 1$ is equal to
$$c_3(q_0+aq_{2k+1}+\frac{a^2}{2}q_{2(2k+1)}+\dots,
q_2+aq_{2+2k+1}+\frac{a^2}{2}q_{2+2(2k+1)}+\dots, \dots).$$
Corollary is proved.
\end{proof}

In the rest of this section we prove that $G_-$ acts transitively on $GW_0$.

\begin{lem}\label{N=1}
Let $F\in GW_0$ and let $c_3=-\frac{1}{6}(q_0+\al_2a_2+\dots)$ be the
coefficient of $q_1^{-1}$ in $F$.
Then
\begin{itemize}
\item  $\al_2=0$;
\item  if $\al_i=0$ for $i\le 2N-1$, $N=1,2,\dots$, then $\al_{2N}=0$.
\end{itemize}
\end{lem}
\begin{proof}
We apply derivatives $\pa_0\pa_2$ to both sides of $(SE(0))$ and obtain
$$\pa_2\pa_0^2c_3+\pa_2\pa_0(q_2c_1)+\sum_{i\ge 2} q_{i+1}\pa_2\pa_0\pa_ic_2=0.$$
This gives $\al_2=0$, because $\pa_0c_1=0$ and $c_2$ is a degree $2$ series.

Now suppose $\al_2=\al_3=\dots=\al_{2N-1}=0$. Let
$$c_2=\frac{1}{2}\sum_{i,j\ne 1} \be_{ij} q_iq_j.$$
The equation $(SE(1))$ gives
\beq\label{c0c1}
c_0=-\be_{2,0},\ c_1=-\be_{3,0}q_2-\be_{4,0}q_3-\dots.
\eeq
Thus, the equation $(SE(0))$ reads as
$$
-\frac{1}{2}(q_0+\al_2q_2+\dots)^2-q_2(\be_{3,0}q_2+\be_{4,0}q_3+\dots)+
\sum_{i\ge 2}q_{i+1}\sum_{j\ge 0} q_j\be_{i,j}=-\frac{1}{2}q_0^2.
$$
Extracting the coefficient of $q_nq_0$ we obtain $\al_2=0$ and
\beq\label{alnbe0}
\al_n=\be_{n-1,0}, n\ge 3.
\eeq
The coefficient of $q_2q_n$, $n\ge 2$ gives
\beq\label{be0be2}
\be_{n+1,0}=\be_{n-1,2}.
\eeq
Finally the coefficient of $q_iq_j$, $i,j\ge 3$ yields
\beq\label{bebe}
\be_{i-1,j}+\be_{i,j-1}=\al_i\al_j.
\eeq
Thus we obtain the following equations:
\begin{multline*}
\al_{2n}=\be_{2n-1,0}=\be_{2n-3,2}=\al_{2n-3}\al_3-\be_{2n-4,3}=\dots\\
=\al_{2n-3}\al_3-\al_{2n-4}\al_4+\dots\pm \frac{1}{2}\al_{n+1}^2.
\end{multline*}
This proves our Lemma.
\end{proof}

\begin{lem}
\label{c3}
A point $F\in GW_0$ is uniquely determined by the $q_1^{-1}$-coefficient $c_3$.
\end{lem}
\begin{proof}
Formulas $(\ref{alnbe0})$, $(\ref{be0be2})$ and $(\ref{bebe})$ show that
$\be_{i,j}$ can be expressed via $\al_i$.
Therefore $c_0,c_1,c_2$ are determined by $c_3$.

We now show that the same is true for $c_n$ with $n\ge 4$.
In fact, the equation $TRR_{0,1,1}(-n-1)$ allows to express $c_n$ in terms
of $\pa_0c_n$
and $c_m$, $m<n$. But using $SE(-n+3)$ we can write $\pa_0 c_n$ via
$c_{n-1}$ and $c_{n-2}$. Lemma is proved.
\end{proof}

\begin{cor}
The group $G_-$ acts transitively on $GW_0$.
\end{cor}
\begin{proof}
Note that the coefficient of $q_1^{-1}$ in $F_0^{pt}$ is equal to $-\frac{q_0^3}{6}$.
Therefore it suffices to show that for any $F\in GW_0$ there exist complex numbers
$a_3,a_5,\dots$ such that the coefficient of $q_1^{-1}$ in
$$\exp(a_3z^{-3}+a_5z^{-5}+\dots)F$$
is equal to  $-\frac{q_0^3}{6}$. We find $a_{2k+1}$ by induction on $k$.
Let $$-\frac{1}{6}(q_0+\al_3q_3+\dots)^3$$
be a coefficient of $q_1^{-1}$ in $F$. Then the coefficient of $q_1^{-1}$ in
$\exp(\al_3z^{-3})F$ is a cube of a linear form that contains no $q_3$ term and
therefore (because of the Lemma $\ref{N=1}$) no $q_4$ term. Now suppose
that the coefficient of $q_1^{-1}$ in $F$ is given by
$$-\frac{1}{6} (q_0+\al_{2k+1}q_{2k+1}+\al_{2k+2}q_{2k+2}+\dots)^3.$$
Then the coefficient of $q_1^{-1}$ in
$\exp(\al_{2k+1}z^{-2k-1})F$ is a cube of a linear form that contains no $q_i$ terms,
$i\le 2k+1$ and therefore (because of the Lemma $\ref{N=1}$) no $q_{2k+2}$ term.
This yields our Corollary.
\end{proof}

\begin{theo}
The set of solutions of $(DE)$, $(SE)$ and $(TRR)$ is isomorphic to $\C^\infty$ via the
map
$$(a_3,a_5,\dots)\mapsto  \exp(\sum_{i\ge 1} a_{2i+1} z^{-2i-1})\cdot F^{pt}_0.$$
\end{theo}

\begin{cor}
The system $(DE)$ + $(SE)$ + $(TRR)$ is equivalent to $(DE)$ + $(SE)$ + $(TRR_{0,1,1})$.
\end{cor}
\begin{proof}
Follows from the proof.
\end{proof}

\begin{rem}
For general $N$ one needs both groups $G_\pm$ to generate the space $GW_0$ starting from
the potential of $N$ points. But for $N=1$ the Lie algebra $\g_+$ (see $(\ref{Ra})$) acts trivially
on $F^{pt}_0$ due to the Faber-Pandharipande relations (see \cite{Gi2}, \cite{FP}).
\end{rem}

\newcounter{a}
\setcounter{a}{1}

\end{document}